\documentclass{mcom-l}

\newtheorem{theorem}{Theorem}[section]
\newtheorem{lemma}[theorem]{Lemma}

\theoremstyle{definition}

\newtheorem{example}[theorem]{Example}

\theoremstyle{remark}
\newtheorem{remark}[theorem]{Remark}

\numberwithin{equation}{section}

\newcommand{\beq}{\begin{equation}}
\newcommand{\eeq}{\end{equation}}

\def\a{\alpha}
\def\b{\beta}
\def\g{\gamma}
\def\d{\delta}
\def\e{\epsilon}

\def\({\left(}
\def\){\right)}

\begin{document}

\title{On using symmetric polynomials for constructing root finding methods}

\author{Dmitry I. Khomovsky}
\address{Faculty of Physics, Lomonosov Moscow State University,  1-2 Leninskie Gory, 119991 Moscow, Russia}
\email{khomovskij@physics.msu.ru}
\thanks{}


\date{}

\dedicatory{Dedicated to the  memory  of  my school teacher Alexander L. Smirnov.}


\begin{abstract}
We propose an approach to constructing iterative methods for finding polynomial roots simultaneously. One feature of this approach is using the fundamental theorem of symmetric polynomials.  Within this framework, we reconstruct many of the existing root finding methods.
The new results presented in this paper are some modifications of the Durand-Kerner method.
\end{abstract}

\maketitle

\section{Introduction}
Let $f(\textrm{z})$ be a polynomial of degree $n$  with coefficients in $\mathbb{C}$ and let its factorization over the complex numbers be $f(\textrm{z})=\prod_{j=1}^{n}(\textrm{z}-\lambda_j)$, where $\lambda_j$ $(j=1,2,...,n)$ are the roots (zeros) of $f(\textrm{z})$.

Let us consider some known methods for simultaneous approximation of roots. The classical (Weierstrass) Durand-Kerner method \cite{Do,Du,K,We} is related to
\beq\label{1}
\textrm{z}_i^{(k+1)}=\textrm{z}_i^{(k)}-\frac{f(\textrm{z}_i^{(k)})}{\prod\limits_{\substack{j=1\\j\not=i}}^{n}(\textrm{z}_i^{(k)}-\textrm{z}_j^{(k)})} \,\,\,\,\,(i=1,\ldots,n),
\eeq
here $k$ is the iteration number. Further in similar formulas we will use $\textrm{z}_i$ and $\hat{\textrm{z}}_i$ instead $\textrm{z}_i^{(k)}$ and $\textrm{z}_i^{(k+1)}$, respectively.  If the roots $\lambda_i$ $(i=1,2,...,n)$ are distinct and the initial approximations $\textrm{z}_i^{(0)}$ $(i=1,2,...,n)$ are close to them, then the method is of quadratic convergence proven by Dochev \cite{Do}.

The Maehly-Ehrlich-Alberth method \cite{A,E,Ma} with cubic convergence\footnote{Here and further we imply only the case of simple roots and good initial approximations.} deals with
\beq\label{2}
\hat{\textrm{z}}_i=\textrm{z}_i-\left[\frac{f'(\textrm{z}_i)}{f(\textrm{z}_i)}-\sum\limits_{j\not=i}\frac{1}{\textrm{z}_i-\textrm{z}_j}\right]^{-1}\,\,\,\,\,(i=1,\ldots,n).
\eeq
In practice, it is convenient to use a formula which does not contain division by a near-zero value $f(\textrm{z}_i)$, since it may lead to loss of accuracy. So the following formula is used:
\[
\hat{\textrm{z}}_i=\textrm{z}_i-f(\textrm{z}_i) \left[f'(\textrm{z}_i)-f(\textrm{z}_i)\sum\limits_{j\not=i}\frac{1}{\textrm{z}_i-\textrm{z}_j}\right]^{-1}.
\]
There are modifications that significantly improve the  iterative schemes above (see Petcovic and Milovanovic \cite{Mi1,Mi2,PM} and references therein).

The Ostrowski-Gargantini method \cite{Gar,O} having the fourth order of convergence is based on the following iterative formula:
\beq\label{3}
\hat{\textrm{z}}_i=\textrm{z}_i - \left[\left(\frac{f'(\textrm{z}_i)}{f(\textrm{z}_i)}\right)^{2}-\frac{f''(\textrm{z}_i)}{f(\textrm{z}_i)}-\sum\limits_{j\not=i}\frac{1}{(\textrm{z}_i-\textrm{z}_j)^2}\right]^{-1/2}_{*}\,\,\,\,\,(i=1,\ldots,n).
\eeq
The symbol  $*$ denotes that one of the values of the square root  (more appropriate) is chosen. In using such notation we follow \cite{PMS,PS}. A  criterion  for  the  choice  of an  appropriate  value  of  the  square  root is given in \cite{Gar}; we need to choose such a value of the square root so that the following is minimal:
\beq\label{4}
\left|\frac{f'(\textrm{z}_i)}{f(\textrm{z}_i)}- \left[\left(\frac{f'(\textrm{z}_i)}{f(\textrm{z}_i)}\right)^{2}-\frac{f''(\textrm{z}_i)}{f(\textrm{z}_i)}-\sum\limits_{j\not=i}\frac{1}{(\textrm{z}_i-\textrm{z}_j)^2}\right]^{1/2}\right|.
\eeq
Since $(\ref{4})$ contains only the terms which must be calculated in the current iteration step, the direct way of choosing a value of the square root, which implies the minimization of $|f(\hat{\textrm{z}}_i)|$, requires more calculations in a general case.

The generalization of $(\ref{2})$, $(\ref{3})$ was presented in \cite{P1,PMS}. This result is as follows:
\beq\label{6}
\hat{\textrm{z}}_i=\textrm{z}_i - \left[F_m(\textrm{z}_i)-\sum\limits_{j\not=i}\frac{1}{(\textrm{z}_i-\textrm{z}_j)^m}\right]^{-1/m}_{*}\,\,\,\,\,(i=1,\ldots,n),
\eeq
where \beq\label{7}F_m(\textrm{z})=\frac{(-1)^{m-1}}{(m-1)!}\frac{d^{m-1}}{d\textrm{z}^{m-1}}\(\frac{f'(\textrm{z})}{f(\textrm{z})}\)\,\,\, (m\in \mathbb{Z}^{+}).\eeq
To choose an appropriate  value  of  the  $m$th  root we can use the minimization of
\beq\label{8}
\left|\frac{f'(\textrm{z}_i)}{f(\textrm{z}_i)}- \left[F_m(\textrm{z}_i)-\sum\limits_{j\not=i}\frac{1}{(\textrm{z}_i-\textrm{z}_j)^m}\right]^{1/m}\right|.
\eeq
The generalized  iterative formula $(\ref{6})$ is locally of $(m+2)$th order of convergence.
For more  information  about simultaneous root-finding methods see \cite{B1,B2,Mc,P2,S}.

In this article, we discuss a new view on iteration methods for the simultaneous
approximation of polynomial roots based on relations for symmetric multivariate polynomials. In the next section, we present a framework to reconstruct all sorts of iterative methods illustrated by some well-known earlier results.
\section{Constructing iterative formulas}
The elementary symmetric polynomials are defined as follows
\beq\label{8.1}
e_0(x_1,\ldots,x_n)=1,\,\,\,\,e_k(x_1,\ldots,x_n)=\sum\limits_{1\leq j_1\dots<j_k\leq n}x_{j_1}\cdots x_{j_k} \,\,\,\,\,(1 \leq k \leq n).
\eeq
It is known that any symmetric polynomial in $x_1,\ldots,x_n$ can be expressed as a polynomial in $e_k(x_1,\ldots,x_n)\,(1\leq k\leq n)$, moreover, such a representation is unique.
\begin{example} For example, we consider the $m$th power sum of $n$ variables, i.e., $p_m(x_1,\ldots,x_n)=\sum_{j=1}^{n}x_j^m$. There is the following recursive procedure:
\begin{align}\label{10}
\nonumber p_1 & = e_1, \\
\nonumber p_2& = e_1 p_1-2e_2, \\
\nonumber p_3&=e_1p_2-e_2p_1+3e_3,\\
\nonumber p_4&=e_1p_3-e_2p_2+e_3p_1-4e_4,\,\,\text{and so on.}
\end{align}
The recurrence relation is
\[
p_m=\sum_{j=1}^{m-1}(-1)^{m-1+j}e_{m-j}p_j+(-1)^{m-1}m e_m,\,\,m\geq 1.
\]
Therefore, we can obtain the representation of $p_m$ via $e_k$ $(1\leq k\leq m)$.
Also, there are explicit formulas which express power sums in terms of elementary symmetric polynomials, see \cite{Me}.
\end{example}
\begin{lemma}
Let $f(\textrm{z})$ be a polynomial of degree $n$  with coefficients in $\mathbb{C}$ and $\lambda_j$ $(j=1,2,...,n)$ be its roots.
For an integer $0\leq k\leq n$ the following holds:
\beq\label{11}
\frac{1}{k!}\frac{f^{(k)}(\textrm{z})}{f(\textrm{z})}=e_k\left(\frac{1}{\textrm{z}-\lambda_{1}},\ldots,\frac{1}{\textrm{z}-\lambda_{n}}\right),
\eeq
here $e_k$ is the elementary symmetric polynomial of degree $k$ in $n$ variables.
\end{lemma}
\begin{proof}
We have the following two formulas which derived from the definition of elementary symmetric polynomials $(\ref{8.1})$:
\beq\label{11.1}
e_k\left(\frac{1}{\textrm{z}-\lambda_{1}},\ldots,\frac{1}{\textrm{z}-\lambda_{n}}\right)\, f(\textrm{z})=
e_{n-k}\left(\textrm{z}-\lambda_{1},\ldots,\textrm{z}-\lambda_{n}\right),
\eeq
\beq\label{11.2}
\frac{d}{d\textrm{z}}\,e_i\left(\textrm{z}-\lambda_{1},\ldots,\textrm{z}-\lambda_{n}\right)=
(n-i+1)\, e_{i-1}\left(\textrm{z}-\lambda_{1},\ldots,\textrm{z}-\lambda_{n}\right).
\eeq
Suppose that $(\ref{11})$ holds for $k=m$ and $m< n$. By using $(\ref{11})$ and $(\ref{11.1})$ we get the following:
\beq\label{11.3}
f^{(m)}(\textrm{z})=m!\,e_{n-m}\left(\textrm{z}-\lambda_{1},\ldots,\textrm{z}-\lambda_{n}\right).
\eeq
From this formula with the help of $(\ref{11.2})$ we obtain
\beq\label{11.4}
f^{(m+1)}(\textrm{z})=(m+1)!\,e_{n-m-1}\left(\textrm{z}-\lambda_{1},\ldots,\textrm{z}-\lambda_{n}\right).
\eeq
Thus, we conclude that $(\ref{11})$ also holds for  $k=m+1$.
For $ k = 0 $ the statement of the lemma is true. Then, using mathematical induction, we complete the proof.
\end{proof}
This lemma is used to construct iterative formulas. The main idea is as follows: suppose we take some symmetric polynomial in the variables $1/(z-\lambda_j)$ $(1\leq j\leq n)$ and express it via elementary symmetric polynomials, then using $(\ref{11})$, we obtain a formula which, after simple transformations, will give us a simultaneous root-finding method.
\begin{example}
Let us consider the polynomial $p_3\left((z-\lambda_1)^{-1},\ldots,(z-\lambda_n)^{-1}\right)$. There is the representation $p_3=e_1^3 - 3 e_2 e_1 + 3 e_3$. Using $(\ref{11})$, we obtain
\[
\sum\limits_{j=1}^{n}\frac{1}{(\textrm{z}-\lambda_j)^3}=\left(\frac{f'(z)}{f(z)}\right)^3-\frac{3 f'(z)f''(z)}{2 f(z)^2}+\frac{f'''(z)}{2f(z)}.
\]
Making simple transformations, we derive an explicit expression for $\lambda_i$. We get
\[
\lambda_i=\textrm{z}-\left[\left(\frac{f'(z)}{f(z)}\right)^3-\frac{3 f'(z)f''(z)}{2 f(z)^2}+\frac{f'''(z)}{2f(z)}-\sum\limits_{j\not=i}\frac{1}{(\textrm{z}-\lambda_j)^3}\right]^{-1/3}.
\]
Finally, we have the following iterative method
\beq\label{14}
\hat{\textrm{z}}_i=\textrm{z}_i-\left[\left(\frac{f'(\textrm{z}_i)}{f(\textrm{z}_i)}\right)^3-\frac{3 f'(\textrm{z}_i)f''(\textrm{z}_i)}{2 f(\textrm{z}_i)^2}+\frac{f'''(\textrm{z}_i)}{2f(\textrm{z}_i)}-\sum\limits_{j\not=i}\frac{1}{(\textrm{z}_i-\textrm{z}_j)^3}\right]^{-1/3}_{*}.
\eeq
This is exactly $(\ref{6})$ when $m=3$. Usually, the method $(\ref{14})$ is not used in practice.
\end{example}
\begin{remark}
If we consider $p_m\left((\textrm{z}-\lambda_1)^{-1},\ldots,(\textrm{z}-\lambda_n)^{-1}\right)$, then we obtain $(\ref{6})$ and derive the following relation
\beq\label{15}F_m(\textrm{z})=u_m\left(\frac{f'(\textrm{z})}{f(\textrm{z})},\frac{1}{2!}\frac{f''(\textrm{z})}{f(\textrm{z})},\ldots,\frac{1}{n!}\frac{f^{(n)}(\textrm{z})}{f(\textrm{z})}\right),
\eeq
where the polynomial $u_m$ is defined by
\beq\label{16}p_m(x_1,\ldots,x_n)=u_m\left(e_1(x_1,\ldots,x_n),\ldots,e_n(x_1,\ldots,x_n)\right).
\eeq
\end{remark}
\subsection*{Halley's method for simultaneous approximation of polynomial zeros}
Let $\a$ and $\b$ be nonzero elements in $\mathbb{C}$.  We consider the symmetric polynomial
$\a p_2+\b p_1^2$ in variables  $1/(z-\lambda_j)$ $(1\leq j\leq n)$. Let us introduce the notation:
\beq\label{17} q = \frac{1}{z-\lambda_i},\,\,\,\,\,\, S_r = \sum\limits_{j\not=i}\frac{1}{\left(\textrm{z}-\lambda_j\right)^r}\,\,(r\in\mathbb{Z}^+).\eeq
Then  $p_r=q^r+S_r$. Using this, we have the following:
\begin{align}\label{18}
\nonumber \a p_2+\b p_1^2 & = \a(q^2+S_2)+\b(q^2+2q S_1+S_1^2) \\
\nonumber & = \a(q^2+S_2)+\b(q^2+2 q (p_1-q)+S_1^2)\\
\nonumber & =(\a-\b)q^2+2\b p_1 q+\a S_2+\b S_1^2.
\end{align}
We see that it is convenient to put $\a=\b=1$. Then $q=(p_2-S_2+p_1^2-S_1^2)/(2 p_1)$.
Since $p_1=e_1$, $p_2=e_1^2-2e_2$, with the help of $(\ref{11})$ we get
\beq\label{19}
\lambda_i=\textrm{z}-\frac{2f(\textrm{z})f'(\textrm{z})}{2{[f'(\textrm{z})]}^2-f(\textrm{z})f''(\textrm{z})-{[f(\textrm{z})]}^2(S_2+S_1^2)}.
\eeq
Finally, this formula leads to the simultaneous root-finding method
\beq\label{20}
\hat{\textrm{z}}_i=
\textrm{z}_i-\frac{2f(\textrm{z}_i)f'(\textrm{z}_i)}{2{[f'(\textrm{z}_i)]}^2-
f(\textrm{z}_i)f''(\textrm{z}_i)-
{[f(\textrm{z}_i)]}^2\bigm(\sum\limits_{j\not=i}(\textrm{z}_i-\textrm{z}_j)^{-2}+
{[{\sum\limits_{j\not=i}(\textrm{z}_i-\textrm{z}_j)^{-1}}]}^2\bigm)}.
\eeq
This result was derived by Wang and Zheng in \cite{Wa}. Since $(\ref{19})$ is related to  Halley's method \cite{Gan} for solving a nonlinear equation, so $(\ref{20})$ is sometimes called the Halley-like method for simultaneous approximation of polynomial zeros. Its convergence analysis can be found in \cite{C,P2}, the method is locally of the fourth order of convergence. In the next section we will get this result.
\begin{remark}
It is necessary to clarify how we came to the idea of choosing the polynomial $\a p_2+\b p_1^2$.
First, we considered the cases when the starting polynomials are $p_2$, $p_1^2$.
In both cases, we obtained fourth-order methods, but they contained a squaring operation, see the Ostrowski-Gargantini method $(\ref{3})$. Then we chose the starting polynomial as a linear combination of $p_2$ and $p_1^2$ in order to exclude a squaring operation by choosing values of the coefficients $\a, \b$. If we deal with $\a, \b$ as symbolic parameters (without setting them equal to certain values), we would get a family of fourth-order methods.
\end{remark}
\subsection*{Simultaneous Householder's method}
In 1984, Wang and Zheng \cite{Wa} presented a family of iterative methods. This family contains the Maehly-Ehrlich-Alberth method $(\ref{2})$, and the Halley-like method $(\ref{20})$; the authors used a concept based on Bell's polynomials. Below, within the proposed framework, we reproduce some results.

We consider $\a p_3+\b p_1 p_2+\g p_1^3$ in variables $1/(z-\lambda_j)$ $(1\leq j\leq n)$. Then
\begin{align}\label{21}
\nonumber  &\a p_3+\b p_1 p_2+\g p_1^3  = \a(q^3+S_3)+\b(q+S_1)(q^2+S_2)+\g (q+S_1)^3\\
\nonumber  &= (\a+\b+\g) q^3+(\b+3\g) S_1 q ^2+(\b S_2+3 \g S_1^2)q+\a S_3+\b S_1S_2+\g S_1^3\\
\nonumber  &=(\a-\b+\g)q ^3+ (\b-3 \g) p_1q ^2+(\b p_2+3 \g p_1^2)q+\a S_3+\b S_1S_2+\g S_1^3.
\end{align}
We put $\a=2, \b=3, \g=1$ and get
\[q=\frac{2 (p_3-S_3)+3 (p_1 p_2-S_1 S_2)+ p_1^3-S_1^3}{3 (p_2+p_1^2)}.\]
Therefore, we have
\beq\label{22}
\lambda_i=\textrm{z}-\frac{6f {f^{\prime}}^2-3f^2 f^{\prime\prime}}{6{f^{\prime}}^3 -6 f f^{\prime}f^{\prime\prime} + f^2f^{\prime\prime\prime}-f^3(2S_3+3S_1S_2+S_1^3)}.
\eeq
Using this formula, we can get the corresponding simultaneous root-finding method, which is connected to Householder's method \cite{H} for solving a nonlinear equation $g(x)=0$, where $g$ is a function in one real variable. Indeed, the iterative formula of the $d$th-order Householder's method\footnote{The rate of convergence of the method has order $d+1.$} is
\beq\label{23}
\hat{x} = x + d\; \frac { \left(1/g\right)^{(d-1)} (x) } { \left(1/g\right)^{(d)} (x) }\,\,\,(d\in\mathbb{Z}^+),
\eeq
then for $d=3$ we have
\[
\hat{x}=x-\frac{6g {g^{\prime}}^2-3g^2 g^{\prime\prime}}{6{g^{\prime}}^3 -6 g g^{\prime}g^{\prime\prime} + g^2g^{\prime\prime\prime}}.
\]

Let us consider $\a p_4+\b p_1 p_3+\g p_2^2+\d p_1^2p_2+\e p_1^4$; the number of summands  is equal to the integer partition of $4.$ By analogy with the previous we get

\smallskip
\noindent$\a(p_4-S_4)+\b(p_1p_3-S_1S_3)+\g(p_2^2-S_2^2)+\d(p_1^2p_2-S_1^2S_2)+\e(p_1^4-S_1^4)$\hfil

\smallskip
\par \hfil \noindent$= (\a-\b-\g+\d-\e) q ^4+(\b-2\d+4\e) p_1 q ^3+\((2\g-\d)p_2+(\d-6\e)p_1^2\)q^2$\hfil

\smallskip
\par \hfill \noindent$+(\b p_3+2 \d p_1p_2+4\e p_1^3)q.$

\smallskip
\noindent We put $\e=1$, then in order to obtain a linear equation with respect to the variable $q$ we need to solve the following system:
\beq\begin{cases}\nonumber
\a-\b-\g+\d-1=0,\\
\b-2\d+4=0,\\
2\g-\d=0,\\
\d-6=0.
\end{cases}\eeq
The solution is $\a=6, \b=8, \g=3, \d=6.$ Finally, we have
\beq\label{24}
\lambda_i=\textrm{z}-\frac{4f( 6{f^{\prime}}^3-6f f^{\prime}f^{\prime\prime}+f^2f^{(3)})}{24{f^{\prime}}^4 -
36 f {f^{\prime}}^2f^{\prime\prime}
+6f^2{f^{\prime\prime}}^2+8f^2f^{\prime}f^{(3)}-f^3f^{(4)}-f^4 T},
\eeq
where $T=6S_4+8S_1S_3+3S_2^2+6S_1^2S_2+S_1^4.$ Since $(\ref{24})$ is also related to $(\ref{23})$, we can represent $(\ref{19}), (\ref{22}), (\ref{24})$ in the following form
\beq\label{25}
\lambda_i=\textrm{z}+ d\; \frac { \left(1/f\right)^{(d-1)} (\textrm{z}) } { \left(1/f\right)^{(d)} (\textrm{z})+(-1)^{d-1}H_d/f(\textrm{z}) },
\eeq
where $d=2,3,4$, respectively and
\begin{align}
\nonumber  H_2&=S_2+S_1^2,\\
\nonumber  H_3&=2S_3+3S_1S_2+S_1^3, \\
\nonumber  H_4&=6S_4+8S_1S_3+3S_2^2+6S_1^2S_2+S_1^4.
\end{align}
Since the relation $(\ref{25})$ is already established \cite{Wa,P2} for any positive integer $d$,
we will not do it in this paper. It should be noted that when $d=1$, we have $H_1=S_1$. This case corresponds  to the Maehly-Ehrlich-Alberth method $(\ref{2})$.
\subsubsection*{The explicit formula for $H_d$}
The homogeneous symmetric polynomial of degree $k$ in $x_1,\ldots,x_n$ is
\[
h_k (x_1, \dots,x_n) = \sum_{1 \leq j_1 \leq \cdots \leq j_k \leq n} x_{j_1} \cdots x_{j_k}.
\]
As is known,  $h_k$ can be expressed in terms of power sums; the formula is as follows:
\beq\label{26}
h_k = \sum\limits_{\substack{r_1 + 2r_2 + \cdots + kr_k = k \\ r_1\ge 0, \ldots, r_k\ge 0}}\,\, \prod_{j=1}^k \frac{p_j^{r_j}}{r_j !\, j^{r_j}}.
\eeq
Using this, we get
\begin{align}
\nonumber  h_2&= (  p_2 + p_1^2 )/2,\\
\nonumber  h_3&= ( 2 p_3 + 3 p_1 p_2 + p_1^3 )/6, \\
\nonumber  h_4&= ( 6 p_4 + 8 p_1 p_3 + 3 p_2^2 + 6 p_1^2 p_2 + p_1^4 )/24.
\end{align}
Therefore, we see that
\begin{align}\label{27}H_d&=d!\,h_d\left(\frac{1}{z-\lambda_1},\ldots,\frac{1}{z-\lambda_{i-1}},\frac{1}{z-\lambda_{i+1}},\ldots,\frac{1}{z-\lambda_n}\right)\\
\nonumber &=\sum\limits_{\substack{r_1 + 2r_2 + \cdots + dr_d = d \\ r_1\ge 0, \ldots, r_d\ge 0}}\,\, \prod_{j=1}^d \frac{d!\,S_j^{r_j}}{r_j !\, j^{r_j}}.
\end{align}
Also, $H_d$ can be represented in terms of the Bell polynomials, see \cite{P2,Wa}.
The simultaneous root-finding method based on $(\ref{25})$, $(\ref{27})$ is:
\beq\label{25.1}
\hat{\textrm{z}}_i=\textrm{z}_i+ d\; \frac { \left(1/f\right)^{(d-1)} (\textrm{z}_i) } { \left(1/f\right)^{(d)} (\textrm{z}_i)+(-1)^{d-1}\hat{H}_{d;i}/f(\textrm{z}_i) } \,\,\,(i=1,\ldots,n),
\eeq
where
\beq\label{27.1}\hat{H}_{d;i}=
d!\,h_d\left(\frac{1}{\textrm{z}_i-\textrm{z}_{1}},\ldots,\frac{1}{\textrm{z}_i-\textrm{z}_{i-1}},\frac{1}{\textrm{z}_i-\textrm{z}_{i+1}},\ldots,\frac{1}{\textrm{z}_i-\textrm{z}_n}\right).
\eeq
The order of convergence of the method is $d+2$.
\section{Convergence analysis}
In this section, we show that simultaneous Halley's method is of the fourth order of convergence. We consider only the case when all the roots of $f(\textrm{z})$ are distinct, in other words, we assume that there exists a positive real number $M$ such that $|\lambda_l-\lambda_k|>M$ for any $l\neq k$.
Let us denote the right side of the formula $(\ref{20})$ by $\varphi_i(\textrm{z}_1,\ldots,\textrm{z}_n)$. To study the convergence of the method we put $\textrm{z}_k=\lambda_k+\a_k \varepsilon$ $(1\leq k\leq n)$, where $\varepsilon$ is real and $\a_k,\ldots,\a_n$ are arbitrary complex numbers. Then we consider the expression $\varphi_i(\lambda_1+\a_1\varepsilon,\ldots,\lambda_n+\a_n\varepsilon)$ as a function of the variable $\varepsilon$.
We note that $(\ref{20})$ is obtained from the exact formula $(\ref{19})$. So if $\textrm{z}_i$ is arbitrary   and the remaining $\textrm{z}_j$ are equal to $\lambda_j$, then $\varphi_i(\textrm{z}_1,\ldots,\textrm{z}_n)=\lambda_i$. In this case only one iterative step is necessary to obtain $\lambda_i$. Thus, a computational error in some iteration step is caused by errors related to the sums in $\varphi_i$. Therefore, it is convenient  to get the following:

\smallskip
\noindent${f(\textrm{z}_i)}^2\bigg(\sum\limits_{j\not=i}{(\textrm{z}_i-\textrm{z}_j)^{-2}}+
\bigm[{\sum\limits_{j\not=i}{(\textrm{z}_i-\textrm{z}_j)^{-1}}}\bigm]^2\bigg)$\hfil

\smallskip
\par \hfill \noindent$= {f(\textrm{z}_i)}^2\bigg(\sum\limits_{j\not=i}{(\textrm{z}_i-\lambda_j)^{-2}}+
\bigm[{\sum\limits_{j\not=i}{(\textrm{z}_i-\lambda_j)^{-1}}}\bigm]^2\bigg)+ O(\varepsilon^3) \text{ as } \varepsilon\to 0.$

\smallskip
\noindent Here, we use that $f(\textrm{z}_i)=f(\lambda_i+\alpha_i\varepsilon)=O(\varepsilon)$ and (since the roots are distinct)
\begin{align}
\nonumber  \sum\limits_{j\not=i}{(\textrm{z}_i-\textrm{z}_j)^{-2}}+
\bigm[\sum\limits_{j\not=i}{(\textrm{z}_i-\textrm{z}_j)^{-1}}\bigm]^2 & = \sum\limits_{j\not=i}{(\textrm{z}_i-\lambda_j)^{-2}}+
\bigm[\sum\limits_{j\not=i}{(\textrm{z}_i-\lambda_j)^{-1}}\bigm]^2+O(\varepsilon).
\end{align}
Also, since the roots are distinct, it follows that $f'(\lambda_i)\neq0$.
Then using this, we obtain
\begin{align}\label{100}
\nonumber  &\frac{2f(\textrm{z}_i)f'(\textrm{z}_i)}{2{[f'(\textrm{z}_i)]}^2-
f(\textrm{z}_i)f''(\textrm{z}_i)-
{[f(\textrm{z}_i)]}^2\bigm(\sum\limits_{j\not=i}(\textrm{z}_i-\lambda_j)^{-2}+
{[{\sum\limits_{j\not=i}(\textrm{z}_i-\lambda_j)^{-1}}]}^2\bigm)+O(\varepsilon^3)}\\
\nonumber  =&\frac{2f(\textrm{z}_i)f'(\textrm{z}_i)}{2{[f'(\textrm{z}_i)]}^2-
f(\textrm{z}_i)f''(\textrm{z}_i)-
{[f(\textrm{z}_i)]}^2\bigm(\sum\limits_{j\not=i}(\textrm{z}_i-\lambda_j)^{-2}+
{[{\sum\limits_{j\not=i}(\textrm{z}_i-\lambda_j)^{-1}}]}^2\bigm)}+O(\varepsilon^4).
\end{align}
Finally, we have $\varphi_i(\lambda_1+\a_1\varepsilon,\ldots,\lambda_n+\a_n\varepsilon)=\varphi_i(\lambda_1,\ldots,\lambda_i+\alpha_i\varepsilon\ldots,\lambda_n)+O(\varepsilon^4)$ as $\varepsilon\to 0$. As discussed above, $\varphi_i(\lambda_1,\ldots,\lambda_i+\alpha_i\varepsilon\ldots,\lambda_n)=\lambda_i$. So we conclude that the method has the fourth order of convergence.

The above illustrates how we can analyze the convergence of iterative methods like $(\ref{20})$. But since we did not give a convergence theorem with error estimates, we refer the reader to \cite{C,PH}.
\section{Modifications of the Durand-Kerner method}
Let $f(\textrm{z})=\textrm{z}^n+a_{n-1}\textrm{z}^{n-1}+\cdots+a_1 \textrm{z}+a_0$ and, as before, $f(\textrm{z})=\prod_{j=1}^{n}(\textrm{z}-\lambda_j)$.
We begin by introducing the following notation:
\beq\label{28}
e_{k;i}=e_k\left(\frac{1}{\textrm{z}-\lambda_1},\ldots,\frac{1}{\textrm{z}-\lambda_{i-1}},\frac{1}{\textrm{z}-\lambda_{i+1}},\ldots,\frac{1}{\textrm{z}-\lambda_n}\right)\,\,\,(0\leq k\leq n-1).
\eeq
Also, for convenience, we assume that if $k\geq n$, then $e_{k;i}=0$. It is easy to see that for $k\geq 1$ the following identity holds:
\beq\label{29}
e_{k}=q e_{k-1;i}+e_{k;i}.
\eeq
Here, as above $q=1/(\textrm{z}-\lambda_i)$; also $e_k$ is the elementary symmetric polynomial in variables $1/(\textrm{z}-\lambda_j)$  $(1\leq j\leq n)$. If we put $k=n$ in $(\ref{29})$, then
\[1/q=\frac{e_{n-1;i}}{e_{n}}.\]
By $(\ref{11})$ and $(\ref{28})$ we have \[e_n=\frac{1}{n!}\frac{f^{(n)}(\textrm{z})}{f(\textrm{z})}=\frac{1}{f(\textrm{z})}\,\,\, \text{and}\,\,\,e_{n-1;i}=\prod_{j\not=i}\frac{1}{\textrm{z}-\lambda_j}.
\]
Finally, we get
\beq\label{30}
\lambda_i=\textrm{z}-f(\textrm{z})/\prod_{j\not=i}(\textrm{z}-\lambda_j).
\eeq
As is seen, this is the main formula for the Durand-Kerner method $(\ref{1})$. Although the derivation of $(\ref{30})$  from the full factorization of $f(\textrm{z})$ is simpler, we have shown the technique that will be used below.

Now we put $k=n-1$ in $(\ref{29})$, then
\beq e_{n-1}=q e_{n-2;i}+e_{n-1;i}.\eeq
Since $f(\textrm{z})e_{n-1}=n \textrm{z}+a_{n-1}$ and $f(\textrm{z})e_{n-1;i}=1/q$, we have the following:
\beq\label{311} n \textrm{z}+a_{n-1}=f(\textrm{z})q e_{n-2;i}+1/q.
\eeq
Dividing this formula by $q$ and taking into account that $1/q = \textrm{z}-\lambda_i$, we obtain
\beq\label{31}
(\textrm{z}-\lambda_i)^2-(n \textrm{z}+a_{n-1})(\textrm{z}-\lambda_i)+f(\textrm{z}) e_{n-2;i}=0.
\eeq
This formula can be used to obtain Weierstrass-like methods. We have two possible ways: the first is to solve the equation $(\ref{31})$ in the variable $\lambda_i$, the second is to use $(\ref{30})$ so that the equation becomes linear, which is to be solved in $\lambda_i$. In addition, we use the following formula, which can be proved by simple transformations,
\beq\label{32}
e_{n-2;i}=(n \textrm{z}-\textrm{z}-\sum_{j\not=i}\lambda_j)/\prod_{j\not=i}(\textrm{z}-\lambda_j).
\eeq
Then, following the second way, we have
\beq\label{33}
\lambda_i=\textrm{z}-\frac{1}{n\textrm{z}+a_{n-1}}\frac{f(\textrm{z})}{\prod_{j\not=i}(\textrm{z}-\lambda_j)}
\left[(n-1)\textrm{z}-\sum_{j\not=i}\lambda_j+\frac{f(\textrm{z})}{\prod_{j\not=i}(\textrm{z}-\lambda_j)}\right].
\eeq
The corresponding iterative method is as follows:
\beq\label{33.1}
\hat{\textrm{z}}_i=\textrm{z}_i-\frac{W_i}{n\textrm{z}_i+a_{n-1}}
\left[n\textrm{z}_i-\sum_{j=1}^{n}\textrm{z}_j+W_i\right],\,\text{where } W_i=\frac{f(\textrm{z}_i)}{\prod_{j\not=i}(\textrm{z}_i-\textrm{z}_j)}.
\eeq
This method is very close to the Durand-Kerner method. The convergence analysis can be performed by analogy with the previous section. If initial approximations are good and all the roots of $f(\textrm{z})$ are distinct, then the method has quadratic convergence. We did some numerical tests to investigate the convergence properties of the new method. Based on the results, it can be said that $(\ref{33.1})$ does not have advantages over $(\ref{1})$.
Nevertheless, we generalize the method obtained. We have
\beq\label{344}e_{n-m}=q e_{n-m-1;i}+e_{n-m;i}.\eeq
The following holds:
\beq\label{34}
  e_{n-m} = \frac{v_m(\textrm{z})}{f(\textrm{z})},\, \text{ where }\, v_m(\textrm{z})=\sum_{l=0}^{m}a_{n-m+l}\binom{n-m+l}{n-m}{\textrm{z}}^{l}
\eeq
and
\beq\label{35}  e_{k;i} =\frac{c_{n-k-1;i}}{\prod_{j\not=i}(\textrm{z}-\lambda_j)},
\eeq
where $c_{m;i}=e_{m}(\textrm{z}-\lambda_1,\ldots,\textrm{z}-\lambda_{i-1},\textrm{z}-\lambda_{i+1},\ldots,\textrm{z}-\lambda_n)$ $(0\leq m \leq n-1)$.
From these formulas and $(\ref{344})$ it follows that
\beq\label{36}
(\textrm{z}-\lambda_{i})^2 c_{m-1;i}-(\textrm{z}-\lambda_{i})v_m(\textrm{z})+\frac{f(\textrm{z})}{\prod_{j\not=i}(\textrm{z}-\lambda_j)}c_{m;i}=0.
\eeq
Using $(\ref{30})$, we obtain a linear equation in $\lambda_{i}$, solving which we find
\beq\label{37}
\lambda_i=\textrm{z}-\frac{1}{v_m(\textrm{z})}\frac{f(\textrm{z})}{\prod_{j\not=i}(\textrm{z}-\lambda_j)}
\left[c_{m;i}+\frac{f(\textrm{z})}{\prod_{j\not=i}(\textrm{z}-\lambda_j)}c_{m-1;i}\right].
\eeq
The first values of $c_{m;i}$ are given below:
\begin{align}
\nonumber  c_{0;i} &= 1, \\
\nonumber  c_{1;i} &= (n-1)\textrm{z}-b_1, \\
\nonumber  c_{2;i} &= (n-1)(n-2)\textrm{z}^2/2-(n-2)b_1\textrm{z}+(b_1^2-b_2)/2,
\end{align}
where $b_k=\sum_{j\not=i}\lambda_j^k$ $(k\in \mathbb{Z}^+)$. The general formula  is
\[
c_{m;i}=\sum_{l=0}^{m}\binom{n-1-m+l}{l}\(\sum\limits_{\substack{r_1 + 2r_2 + \cdots + (m-l)r_{m-l} = m-l \\ r_1\ge 0, \ldots, r_{m-l}\ge 0}}\,\, \prod_{j=1}^{m-l} \frac{(-b_j)^{r_j}}{r_j !\, j^{r_j}}\){\textrm{z}}^{l}.
\]
In this formula, we assume that the sum over $r_1,\ldots,r_{m-l}$ is equal to $1$ if $m-l = 0$.
\section*{Acknowledgments}
The author thanks the referees for their helpful suggestions. Thanks also to Prof. Miodrag S. Petkovi\'{c} for pointing to the reference \cite{Wa} and for valuable comments.

\bibliographystyle{amsplain}

\end{document}